\documentclass[14pt]{article}
\usepackage[cp1251]{inputenc}
\usepackage[english]{babel}
\usepackage{amsfonts,amssymb,amsmath,amstext,amsbsy,amsopn,amscd,amsthm,graphicx,euscript}
\usepackage{graphics}
\usepackage{amsthm}
\newtheorem{thm}{Theorem}

\newcommand{\eps}{\varepsilon}

\newcounter{tdfn}
\newenvironment{dfn}
{\vspace{0.15cm}{\bf Definition \arabic{tdfn}.}} {\par
\addtocounter{tdfn}{1}}
\newcounter{trk}
\newenvironment{rk}
{\vspace{0.15cm}{\bf Remark \arabic{trk}.}} {\par
\addtocounter{trk}{1}}
{\endtrivlist}

\def\:{\colon}

\def\R{{\mathbb R}}

\def\0{{\mathbf 0}}
\def\1{{\mathbf 1}}

\def\R{{\mathbb R}}

\title{Non-Reidemeister Knot Theory and Its Applications in Dynamical Systems, Geometry, and Topology}

\author{Vassily Olegovich Manturov \footnote{The present paper is supported by the Laboratory
of Quantum Topology of the Chelyabinsk State Unversity (Russian
Federation government grant 14.Z50.31.0020), by RF President NSh —
1410.2012.1,
 and
by grants of the Russian Foundation for Basic Resarch,
13-01-00830,14-01-91161, 14-01-31288.}}

\date{}

 \def\R{{\mathbb R}}

\begin{document}

\maketitle

\begin{abstract}

Classical knot theory deals with {\em diagrams} and {\em
invariants}. By means of horizontal {\em trisecants}, we construct a
new theory of classical braids with invariants valued in {\em
pictures}.

These pictures are closely related to diagrams of the initial
object.

The main tool is the notion of {\em free $k$-braid group}. In the
simplest case, for free $2$-braids, the word problem and the
conjugacy problem can be solved by finding the minimal
representative, which can be thought of as a graph, and is unique,
as such.

We prove a general theorem about invariants of dynamical systems
which are valued in such groups and hence, in pictures.

We describe various applications of the above theory: invariants of
weavings (collections of skew lines in $\R^{3}$), and many other
objects in geometry and topology.

In general, provided that for some topological objects (considered
up to isotopy, homotopy etc) some easy axioms (coming from some
dimensional constraints) hold, one can construct similar dynamical
systems and picture-valued invariants.

 {\bf Keywords:} Dynamical System, Braid, Knot, Graph, Trisecant, Group, Picture, Invariant, Reidemeister Move, Weave.
\end{abstract}

{\bf AMS MSC} 05C83, 57M25, 57M27,54H20, 05E18

\section{Introduction}

Usually, invariants of mathematical objects are valued in numerical
or polynomial rings, rings of homology groups, etc. In the present
paper, we prove a general theorem about invariants of {\em dynamical
systems} which are valued in pictures, the so-called {\em free
$k$-braids.}

Formally speaking, free  $k$-braids form a group presented by
generators and relations; it has lots of picture-valued invariants;
these invariants (in our new terminology, {\em $2$-braids}). For
free $2$-braids, the following principle can be realized

{\em if a braid diagram $D$ is complicated enough, then it realizes
itself as a subdiagram of any diagram $D'$ contained in $D$}.

Our invariant is constructed by using horizontal trisecant lines.
Herewith,  the set of {\em cricical values} (corresponding to these
trisecants) leads to  a certain picture which appears in all
diagrams equivalent to the initial picture.

The main theorem of the present paper has various applications in
knot theory, geometry, and topology.

It is based on the following important principle:

{\em if in some science an object has intersections of multiplicity
$k$, then such an object can be studied by $k$-braids and their
generalizations.}

In topology, this principle was first demonstrated in terms of {\em
parity} for the case of {\em virtual knots}, see \cite{Parity}.

Given a topological space $\Sigma$, called the {\em configuration
space}; the elements of $\Sigma$ will be referred to as {\em
particles}. The topology on $\Sigma^{N}$ defines a natural topology
on the space of all continuous mappings $[0,1]\to \Sigma^{N}$; we
shall also study mappings $S^{1}\to \Sigma^{N}$, where
$S^{1}=[0,1]/\{0=1\}$ is the circle; these mappings will naturally
lead to {\em closures} of $k$-braids.

Let us fix positive integers $n$ and $k$.

The   {\em space of admissible dynamical systems} ${\cal D}$ is a
closed subset in the space of all maps $[0,1]\to \Sigma^{n}$.

By a {\em dynamical system} we mean an element of $D$. By a {\em
state} of $D$ we mean the ordered set of particles $D(t)\in
\Sigma^{n}$. Herewith, $D(0)$ and $D(1)$ are called {\em the initial
state} and {\em the terminal state}.

 We shall also deal with {\em cyclic
dynamical system}, where $D(0)=D(1)$.

As usual, we shall fix the initial state and the terminal state and
consider the set of admissible dynamics with such initial and
terminal states.

We say that a property ${\cal P}$ of defined for subsets of the set
of $n$ particles from $\Sigma$ is {\em $k$-good}, if the following
conditions hold:

\begin{enumerate}

\item if this property holds for some set of particles then it holds for every subset of this set;

\item this property holds for every set consisting of $k-1$
particles among $n$ ones (hence, for every smaller set);

\item fix $k$ pairwise distinct numbers $i_{1},\dots,i_{k}, i_{j}=1,\dots,n$;
if the property ${\cal P}$ holds for particles with numbers
 $i_{1},\dots, i_{k-1}$ and for the set of particles with numbers $i_{2},\dots, i_{k}$, then it holds
 for the set of all
$k$ particles $i_{1},\dots, i_{k}$.

\end{enumerate}

\begin{rk} Besides {\em statically good properties} ${\cal P}$,
which are defined for subsets of the set of particles regardless the
state of the dynamics, one can talk about {\em dynamically good}
properties, which can be defined for states considered in time;
however, we shall not need it in the present paper.
\end{rk}

The simplest example is the motion of points in the Euclidian space
$\R^{N}$. Here for $k$-good property we take the property of points
to lie on a fixed $(k-2)$-plane ($k\le N$).

\begin{rk} Our main example deals with the case $k=3$, where for particles we
take different points on the plane, and ${\cal P}$ is the property
of points to belong to the same line. In general, particles may be
more complicated objects than just points.
\end{rk}

Let ${\cal P}$ be a $k$-good property defined on a set of  $n$
particles $n>k$. For each  $t\in [0,1]$ we shall fix the
corresponding state of particles, and pay attention to those $t$ for
which there is a set $k$ particles possessing ${\cal P}$; we shall
refer to these moments as  ${\cal P}$-critical (or just {\em
critical}).

\begin{dfn} We say that a dynamical system $D$ is {\em pleasant}, if the set of its critical moments is finite whereas
each for critical moment there exists exactly one $k$-index set for
which the set ${\cal P}$ holds (thus, for larger sets the property
${\cal P}$ does not hold). Such an unordered $k$-tuple of indices
will be called a {\em multiindex} of critical moments.
\end{dfn}

With a pleasant dynamical system we associate its {\em type}
$\tau(D)$, which will of the set of multiindices $m^{1},\dots,
m^{N}$, written as $t$ increases from $t=0$ to $t=1$.

For each dynamics $D$ and each multiindex $m=(m_{1},\dots, m_{k})$,
let us define the {\em $m$-type} of $D$ as the ordered set $t_{1}<
t_{2}< \dots$ of values of $t$ for which the set of particles
$m_{1},\dots, m_{k}$ possesses the property ${\cal P}$.

Notation: $\tau_{m}(D)$. If the number $l$ of values $t_{1}<\dots<
t_{l}$ is fixed, then the type can be thought of as a point in
$\R^{l}$ with coordinates $t_{1},\dots, t_{l}$.

\begin{dfn}
By the {\em type} of a dynamical system, we mean the set of all its
types $\tau_{m}(D)$. Notation: $\tau(D)$.
\end{dfn}

If $D$ is pleasant, then these set are pairwise disjoint.

\begin{dfn} Fix a number $k$ and a $k$-good property  ${\cal P}$. We say that  $D$
is  {\em ${\cal P}$-stable}, if there is a neighbourhood $U(D)$,
where each dynamical system $D'\in U$ is pleasant, whereas for each
multiindex $m$ (consisting of $k$ indices), the number $l$ of
critical values corresponding to this multiindex $U(D)$, the type
$\tau_{m}$ is a continuous mapping $U(D)\to \R^{l}$.
\end{dfn}

We shall often say {\em pleasant dynamical systems} or {\em stable
dynamical systems} without referring to ${\cal P}$ if it is clear
from the context which ${\cal P}$ we mean.

\begin{dfn} By a
{\em deformation} we mean a continuous path $s:[0,1]\to {\cal D}$ in
the space of admissible dynamical systems, from a stable pleasant
dynamics $s(0)$ to another pleasant stable dynamical system $s(1)$.
\end{dfn}

\begin{dfn}
We say that a deformation $s$ is {\em admissible} if:
\begin{enumerate}

\item the set of values $u$, where $s(u)$ is not pleasant or is not stable, is finite,
and for those $u$ where $s(u)$ is not pleasant, $s(u)$ is stable.

\item inside the stability intervals, the $m$-types are
continuous for each multiindex $m$.

\item for each value $u=u_{0}$, where $s(u_{0})$ is not pleasant, exactly one of the two following
cases occurs:

\begin{enumerate}
\item There exists exactly one $t=t_{0}$ and exactly one $(k+1)$-tuple
$m=(m_{1},\dots, m_{k+1})$ satisfying ${\cal P}$ for this $m$
(hence, ${\cal P}$ does not hold for larger sets).

Let ${\tilde m}_{j}=m\backslash \{m_{j}\},j=1,\dots, k+1$  For types
$\tau_{{\tilde m}_{j}}$, choose those coordinates $\zeta_{j}$, which
correspond to the value $t=t_{0}$. It is required that for all these
values  $u_{0}$ all functions $\zeta_{j}(u)$ are smooth, and all
derivatives $\frac{\partial z_{j}}{\partial{u}}$ are pairwise
distinct;

\item there exists exactly one value $t=t_{0}$ and exactly two
multiindices  $m=\{m_{1},\dots, m_{k}\}$ and $m'=\{m'_{1},\dots,
m'_{k}\}$ for which  ${\cal P}$  holds; we require that $Card(m\cap
m')<k-1$.
\end{enumerate}

\item For each value $u$, where the dynamical system $s(u)$ is not
stable, there exists a value $t=t_{0}$, which is not critical for
$D_{u}$, and a multiindex $\mu=(\mu_{1},\dots,\mu_{k})$, for which
the following holds. For some small  $\eps$ all dynamical systems
$D_{u_0}$ for $u_{0}\in (u-\eps,u)$ and $u_{0}\in (u,u+\eps)$ are
stable (for $\delta<\eps$), and the type $\tau_{\mu}(D_{u+\delta})$
differs from the type $\tau_{\mu}(D_{u-\delta})$ by an
addition/removal of two identical multiindices  $\mu$ in the
position corresponding to $t_{0}$.

\end{enumerate}
\end{dfn}

For the space of deformation, one defines an induced topology.

\begin{dfn} We say that a $k$-good property ${\cal P}$ is $k$-{\em correct} for the space of admissible dynamical
systems, if the following conditions hold:

\begin{enumerate}

\item In each neighbourhood of any dynamical system $D$ there exists a pleasant dynamical system $D'$.

\item For each deformation $s$ there exists an {\em admissible
} deformation with the same ends $s'(0)=s(0),s'(1)=s(1)$.

\end{enumerate}

\end{dfn}

\begin{dfn} We say that two dynamical systems $D_{0},D_{1}$ are {\em equivalent}, if there exists a deformation $s, s(0)=D_{0},s(1)=D_{1}$.
\end{dfn}

Thus, if we talk about a correct  ${\cal P}$-property, we can talk
about an admissible deformation when defining the equivalence.

\section{Free $k$-braids}

Let us now pass to the definition of the {\em $n$-strand free
$k$-braid group $G_{n,k}$}.

Consider the following $\left(\stackrel{n}{k}\right)$ generators
$a_{m}$, where $m$ runs the set of all unordered $k$-tuples
$m_{1},\dots, m_{k},$ whereas each $m_{i}$ are pairwise distinct
numbers from $\left\{1,\dots, n\right\}$.

For each $(k+1)$-tuple $U$ of indices $u_{1},\dots, u_{k+1} \in
\{1,\dots, n\}$,consider the $k+1$ sets $m^{j}=U\backslash
\{u_{j}\}, j=1,\dots, k+1$. With $U$, we associate the relation

$$a_{m^1}\cdot a_{m^2}\cdots a_{m^{k+1}}= a_{m^{k+1}}\cdots a_{m^2}\cdot a_{m^1}; \eqno(1)$$
for two tuples $U$ and ${\bar U}$, which differ by order reversal,
we get the same relation.

Thus, we totally have
$\frac{(k+1)!\left(\stackrel{n}{k+1}\right)}{2}$ relations.

We shall call them the {\em tetrahedron relations}.

For $k$-tuples  $m,m'$ with $Card(m\cap m')<k-1$, consider the {\em
far commutativity relation:}

$$a_{m}a_{m'}=a_{m'}a_{m}\eqno(2).$$

Note that the far commutativity relation can occur only if $n>k+1$.

Besides that, for all multiindices $m$, we write down the following
relation:

$$a_{m}^{2}=1 \eqno(3)$$

Define $G_{n}^{k}$ as the quotient group of the free group generated
by all $a_{m}$ for all multiindices $m$  by relations (1), (2) and
(3).

{\bf Example.} The group $G_{3}^{2}$ is $\langle
a,b,c|a^{2}=b^{2}=(abc)^{2}=1\rangle,$ where
$a=a_{12},b=a_{13},c=a_{23}$.

Indeed, the relation $(abc)^{2}=1$ is equivalent to the relation
$abc=cab=1$ because of $a^{2}=b^{2}=c^{2}=1$. This obviously yields
all the other tetrahedron relations.

{\bf Example.} The group $G_{4}^{3}$ is isomorphic to $\langle
a,b,c,d|a^{2}=b^{2}=c^{2}=1,(abcd)^{2}=1,(acdb)^{2}=1,(adbc)^{2}=1\rangle$.
Here $a=a_{123},b=a_{124},c=a_{134},d=a_{234}.$

It is easy to check that instead of $\frac{4!}{2}=12$ relations, it
suffices to take only $\frac{3!}{2}$ relations.

By the {\em complexity} of a word we mean the number of letters in
this word, by the {\em complexity of a free $k$-braid} we mean the
minimal complexity of all words representing it. Such words will be
called {\em minimal representatives}. The tetrahedron relations (in
the case of free $2$-braids we call them the {\em triangle
relations}) and the far commutativity relations do not change the
complexity, and the relation $a_{m}^{2}=1$ increases or decreases
the complexity by $2$.

As usual in the group theory, it is natural to look for minimal
complexity words represetning the given free $k$-braid.

If we deal with conjugacy classes of free $k$-braids, one deals with
the complexity of cyclic words.

The number of words of fixed complexity in a finite alphabet is
finite; $k$-braids and their conjugacy classes are the main
invariant of the present work.

Let us define the following two types of homomorphisms for free
braids. For each $l=1,\dots, n$, there is an  {\em index forgetting
homomorphism} $f_{l}: G_{n}^{k}\to G_{n-1}^{k-1}$; this homomorphism
takes all generators $a_{m}$ with multiindex $m$ not containing $l$
to the unit element of the group, and takes the other generators
$a_{m}$ to $a_{m'}$, where $m'=m\backslash \{l\}$; this operation is
followed by the index renumbering.

The  {\em strand-deletion homomorphism} $d_{j}$ is defined as a
homomorphism  $G_{n}^{k}\to G_{n-1}^{k}$; it takes all generators
$a_{m}$ having multiindex containing $j$ to the unit element; after
that we renumber indices.

The free $2$-braids (called also {\em pure free braids}) were
studied in
 \cite{Parity,MW, NewBracket}.

For free $2$-braids, the following theorem holds.
\begin{thm}
Let $b'$ be a word representing a free $2$-braid $\beta$. Then every
word $b$ which is a minimal representative of  $\beta$, is
equivalent by the triangle relations and the far commutativity
relation to some subword of the word $b'$.

Every two representatives $b_{1}$ and $b_{2}$ of the same free
$2$-braid  $\beta$ are equivalent by the triangle relations and the
far commutativity relations.\label{theo2}
\end{thm}

Thus, for free $2$-braids, the recognition problem can be solved by
means of considering its minimal representative.

The main idea of the proof of this theorem is similar to the
classification of homotopy classes of curves in $2$-surfaces due to
Hass and Scott \cite{HS}: in order to find a minimal representative,
one looks for ``bigon reductions''  until possible, and the final
result is unique up to third Reidemeister moves for the exception of
certain special cases (multiple curves etc). For free $2$-braids
``bigon reductions'' refer to some cancellations of generators
similar generators $a_{m}$ and $a_{m}$ in good position with respect
to each other, see Fig. \ref{bigon}, third Reidemeister moves
correspond to the triangle relations, and the far commutativity does
not change the picture at all \cite{NewBracket}.

For us, it is crucial to know that {\em when looking at a free
$2$-braid, one can see which pairs of crossings can be cancelled.}
Once we cancel all possible crossings, we get an invariant picture.

Thus, {\em we get a complete picture-valued invariant of a free
$2$-braid}.

Theorem \ref{theo2} means that this picture (complete invariant)
occurs as a sub-picture in every picture representing the same free
$2$-braid.

A complete proof of Theorem \ref{theo2} will be given in a separate
publication.

However, various homomorphisms  $G_{n}^{k}\to G_{n-1}^{k-1}$, whose
combination lead to homomorphisms of type $G_{n}^{k}\to
G_{n-k+2}^{2}$, allow one to construct lots of invariants of groups
$G_{n}^{k}$ valued in {\em picture}.

In particular, these pictures allow one to get easy estimate for the
complexity of braids and corresponding dynamics.

Free $k$-braids will be considered in a separate publication.

\section{The Main Theorem}

Let ${\cal P}$ be a  $k$-correct property on the space of admissible
dynamical systems with fixed initial and final states.

Let  $D$ be a pleasant stable dynamical system decribing the motion
of $n$ particles with respect to ${\cal P}$. Let us enumerate all
critical values $t$ corresponding to all multiindices for $D$, as
$t$ increases from $0$ to $1$. With $D$ we associate an element
$c(D)$ of $G_{n}^{k}$, which is equal to the product of $a_{m}$,
where $m$ are multiindices corresponding to critical values of $D$
as $t$ increases from $0$ to $1$.

\begin{thm}
Let $D_{0}$ and $D_{1}$ be two equivalent stable pleasant dynamics
with respect to ${\cal P}$. Then $c(D_{0})=c(D_{1})$ are equal as
elements of $G_{n}^{k}$.
\end{thm}

\begin{proof}
Let us consider an admissible deformation $D_{s}$ between $D_{0}$
and $D_{1}$. For those intervals of values $s$, where $D_{s}$ is
pleasant and stable, the word representing $c(D_{s})$, does not
change by construction. When passing through those values of $s$,
where $D_{s}$ is not pleasant or is not stable, $c(D(s))$ changes as
follows:

\begin{enumerate}
\item Let $s_{0}$ be the value of the parameter deformation, for
which the property ${\cal P}$ holds for some $(k+1)$-tuple of
indices at some time $t=t_{0}$. Note that $D_{s_0}$ is stable.

Consider the multiindex $m=(m_{1},\dots, m_{k+1})$ for which ${\cal
P}$ holds at $t=t_{0}$ for $s=s_{0}$. Let ${\tilde
m}_{j}=m\backslash m_{j},j=1,\dots, k+1$ For types $\tau_{{\tilde
m}_{j}}$, let us choose those coordinates $\zeta_{j}$, which
correspond to  the intersection $t$ at $s=s_{0}$. As $s$ changes,
these type are continuous functions with respect to $s$.

Then for small $\eps$, for $u=u_{0}+\eps$, the word $c(D_{s})$ will
contain a sequence of letters $a_{{\tilde m}_{j}}$ in a certain
order. For values  $s=s_{0}-\eps$, the word $c(D_{s})$ will contain
the same set of letters in the reverse order. Here we have used the
fact that $D_{s}$ is stable.

\item If for some  $s=s_{0}$ we have a critical value with two different $k$-tuples $m,m'$
possessing ${\cal P}$ and $Card(m\cap m')<k-1$, then the word
$c(D_{s})$ undergoes the relation (2) as $s$ passes through $s_{0}$;
here we also require the stability of $D_{s_0}$.

\item If at some $s=s_{0}$, the deformation $D_{s}$ is unstable, then $c(D_{s})$
changes according to (3) as $s$ passes  through $s_{0}$. Here we use
the fact that the deformation is admissible.

\end{enumerate}
\end{proof}

Let us now pass to our main example, the classical braid group. Here
distinct points on the plane are particles. We can require that
their initial and final poisitions they are uniformely distributed
along the unit circle centered at $0$.

For ${\cal P}$, we take the property to belong to the same line.
This property is, certainly, $3$-good. Every motion of points where
the initial state and the final state are fixed, can be approximated
by a motion where no more than $3$ points belong to the same
straight line at once, and the set of moments where three points
belong to the same line, is finite, moreover, no more than one set
of $3$ points belong to the same line simultaneously. This means
that this dynamical system is pleasant.

Finally, the correctedness of ${\cal P}$ means that if we take two
isotopic braids in general position (in our terminology: two
pleasant dynamical systems connected by a deformation), then by a
small perturbation we can get an admissible deformation for which
the following holds. There are only with finitely many values of the
parameter $s$ with four points on the same line or two triples of
points on the same line at the same moment; moreover, for each such
$s$ only one such case occurs exactly for one value of $t$.

In this example, as well as in the sequel, the properties of being
{\em pleasant} and {\em correct} are based on the fact that every
two {\em general position} state can be connected by a curve passing
through states of codimension $1$ (simplest generation) finitely
many times, and every two paths with fixed endpoints, which can be
connected by a deformation, can be connected by a general position
deformation where points of codimensions $1$ and $2$ occur, the
latter happen only finitely many times.

In particular, the most complicated condition saying that the set of
some $(k+1)$ particles satisfies the property ${\cal P}$, the
corresponding derivatives are all distinct, is also a generial
position argument. For example, assuming that some $4$ points belong
to the same horizontal line (event of codimension $2$), we may
require there is no coincidence of any further parameters (we avoid
events of codimension $3$).

From the definition of the invariant $c$, one easily gets the
following
\begin{thm}
Let $D$ be a dynamical system corresponding to a classical braid.
Then the number of horizontal trisecants of the braid $D$ is not
smaller than the complexity of the free $3$-braid $\beta=c(D)$.
\end{thm}

Analogously, various geometrical properties of dynamical systems can
be analyzed by looking at complexities of corresponding groups of
free $k$-braids, if one can define a $k$-correct property for these
dynamics, which lead to invariants valued in free $k$-braids.

Let us now collect some situations where the above methods can be
applied.

\begin{enumerate}

\item An evident invariant of closed pure classical braids is the conjugacy class of the group $G_{n}^{3}$.
To pass from arbitrary braids to pure braids, one can take some
power of the braid in question.

\item Note that the most important partial case for $k=2$ is the classical Reidemeister braid theory.
Indeed, for a set of points on the plane $Oxy$, we can take for
${\cal P}$ the property that the $y$-coordinates of points coincide.
Then, considering a braid as motion of distinct points in the plane
$z=1-t$ as $t$ changes from $0$ to $1$, we get a set of curves in
space whose projection to $Oxz$ will have intersections exactly in
the case when the property ${\cal P}$ holds. The additional
information coming from the $x$ coordinate, leads one to the
classical braid theory.

\item For classical braids, one can construct invariants for  $k=4$ in a way similar to $k=3$.
In this case, we again take ordered sets of $n$ points on the plane,
and the property ${\cal P}$ means that the set of points belongs to
the same circle or straight line; for three distinct points this
property always holds, and the circle/straight line is unique.

\item With practically no changes this theory can be used for the study of
weavings \cite{Weaves}, collections of projective lines in
$\R{}P^{3}$ considered up to isotopy. Here ${\cal P}$ is the
property of a set of points to belong to the same projective line.
The main difficulty here is that the general position deformation
may contain three lines having infinitely many common horizontal
trisecants. Another difficulty occurs when one of our lines becomes
horizontal; this leads to some additional relations to our groups
$G_{n}^{3}$, which are easy to handle.

\item In the case of points on a $2$-sphere we can define ${\cal P}$ to be the property of
points to belong to the same geodesic. This theory works with an
additional restriction which forbids antipodal points. Some
constraints should be imposed in the case of $2$-dimensional
Riemannian manifolds: for the space of all dynamical systems, we
should impose the restrictions which allow one to detect the
geodesic passing through two points in a way such that if two
geodesics chosen for $a,b$ and for $b,c$ coincide, then the same
geodesic should be chosen for $a,c$.

\item In the case of $n$ non-intersecting projective $m$-dimensional planes in $\R{}P^{m+2}$
considered up to isotopy, the theory works as well. Here, in order
to define the dynamical system, we take a one-parameter family of
projective hyperplanes in general position, for particles we take
$(m-1)$-dimensional planes which appear as intersections of the
initial planes with the hyperplane.

The properties of being good, correct etc. follow from the fact that
in general position ``particles'' have a unique same secant line (in
a way similar to projective lines, one should allow the projective
planes not to be straight).

For $m=1$ one should take $k=3$, in the general case one takes
$k=2m+1$.

This theory will be developed in a separate publication.

\item This theory can be applied to the study of fundamental groups of various discriminant
spaces, if such spaces can be defined by several equalities and
subsets of these equalities can be thought of as property ${\cal
P}$.

\item

The case of classical knots, unlike classical links, is a bit more
complicated: it can be considered as a dynamical system, where the
number of particles is not constant, but it is rather allowed for
two particles to be born from one point or to be mutually
contracted.

The difficulty here is that knots do not possess a group structure,
thus, we don't have a natural order on the set of particles.
Nevertheless, it is possible to construct a map from classical knots
to {\em free $3$-knots} (or free $4$-knots) and study them in a way
similar to free $3$-braids (free $4$-braids).

This theory will be developed in detail in a separate publication.

\end{enumerate}

\begin{rk}
In the case of sets of points in space of dimension $3$, the
property of some points to belong to a $2$-plane (or
higher-dimensional plane) is not correct. Indeed, if three points
belong to the same line, then whatever fourth point we add to them,
the four points will belong to the same plane, thus we can get
various multiindices of $4$ points corresponding to the same moment.

The ``triviality'' of such theory taken without any additional
constraints has the simple description that the configuration space
of sets of points in $\R^{3}$ has trivial fundamental group.
\end{rk}

\section{Pictures}

The free $k$-braids can be depicted by strands connecting points
$(1,0),\dots, (n,0)$ on to points $(1,1),\dots, (n,1)$; every strand
connects $(i,0)$ to $(i,1)$; its projection to the second coordinate
is a homeomorphism. We mark crossings corresponding $a_{i,j,k}$ by a
solid dot where strands $\#i,\#j,\#k$ intersect transversally.

\begin{figure}
\centering\includegraphics{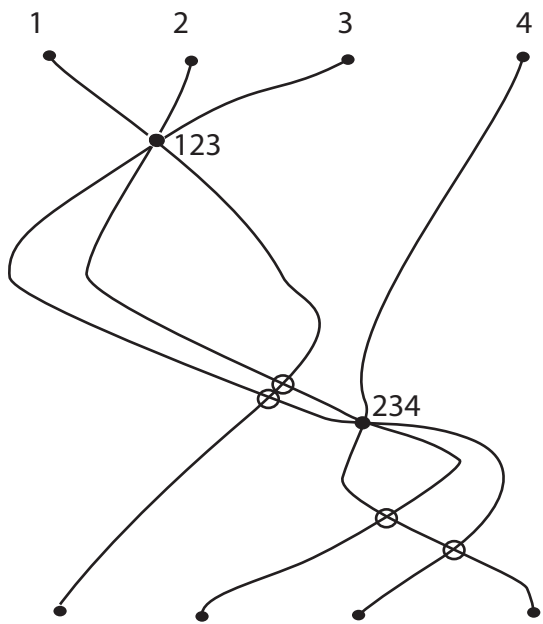} \caption{The $3$-braid
$a_{234}a_{123}$} \label{3braid}
\end{figure}

All other crossings on the plane are artifacts of planar drawing;
they do not correspond to any generator of the group, and they are
encircled. View Fig. \ref{3braid}.

\begin{figure}
\centering\includegraphics{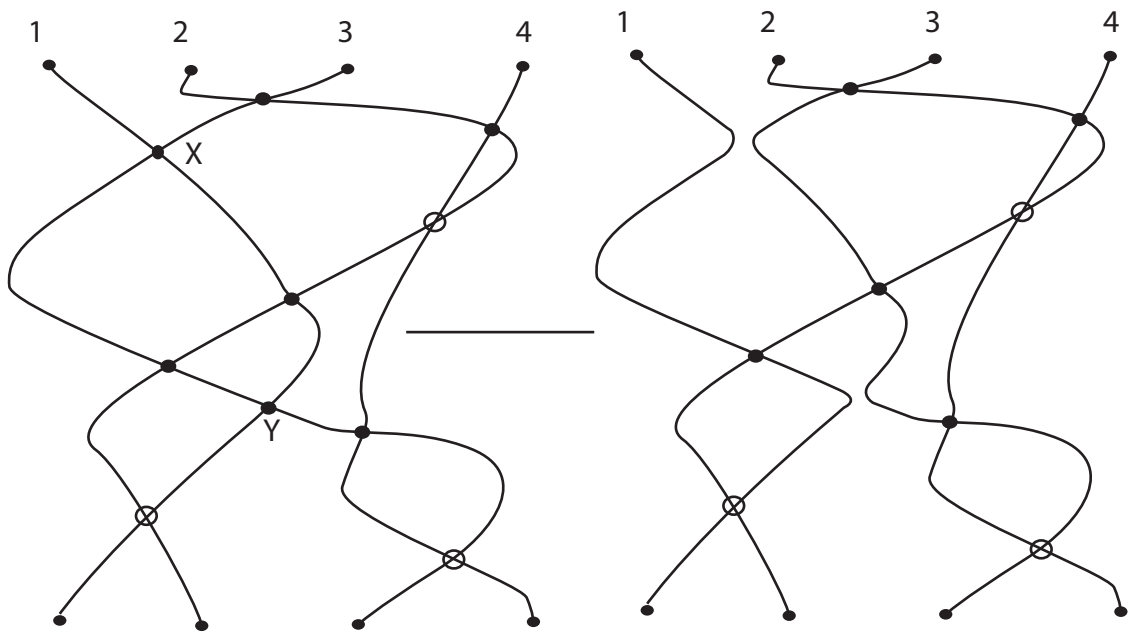} \caption{The bigon
reduction}\label{bigon}
\end{figure}

The clue for the recognition of free $2$-braids is the {\em bigon
reduction} shown in Fig. \ref{bigon}.

Here we reduce the bigon whose vertices vertices are $X,Y$.

The graph which appears after all possible bigon cancellations (with
opposite edge structure at vertices) is the picture which is the
complete invariant of the free $2$-braid (see also
\cite{NewBracket}).

I am grateful to E.I.Manturova, D.P.Ilyutko, I.M.Nikonov, and
V.A.Vassiliev for fruitful discussions.

{\tt

Vassily Olegovich Manturov,

Bauman Moscow State Technical University, Moscow, Russian Federation

vomanturov@yandex.ru

}

\end{document}